\RequirePackage{fix-cm}
\documentclass[twocolumn]{svjour3}          % twocolumn
\smartqed  % flush right qed marks, e.g. at end of proof
\usepackage{graphicx}
\usepackage{amsmath}
\usepackage{amsfonts}
\usepackage{amssymb}
\usepackage{gensymb}
\usepackage{bm}
\usepackage{hyperref}
\usepackage{xcolor} 
\hypersetup{colorlinks, linkcolor={black}, citecolor={black}, urlcolor={black}}
\usepackage{float}
\usepackage{tabularx}
\usepackage{multirow}
\usepackage{array}
\usepackage{makecell}
\usepackage{arydshln}
\usepackage{comment}
\usepackage[numbers]{natbib}
\bibpunct{[}{]}{,}{n}{,}{,}

\newcommand{\abs}[1]{\lvert#1\rvert}
\newcounter{siqiang}
\numberwithin{siqiang}{section}

\begin{document}

\title{A Survey on Machine Learning Solutions for Graph Pattern Extraction%\thanks{Grants or other notes
}


\author
{	
	Kai Siong Yow$^{\ast,\dagger}$\thanks{$^{\ast}$Correspondence:~\texttt{kaisiong.yow@ntu.edu.sg},\\ \texttt{siqiang.luo@ntu.edu.sg}\\ $^{\dagger}$Both authors contributed equally to this research}
	\and
	Ningyi Liao$^{\dagger}$
	\and
	Siqiang Luo$^{\ast}$
	\and
	Reynold Cheng
	\and
	Chenhao Ma
	\and
	Xiaolin Han
}


\institute{Kai Siong Yow, Ningyi Liao and Siqiang Luo \at
              School of Computer Science and Engineering, College of Engineering, Nanyang Technological University, Singapore\\
              (KS Yow is currently on Leave of Absence from the Department of Mathematics and Statistics, Faculty of Science, Universiti Putra Malaysia, 43400 UPM Serdang, Selangor, Malaysia; Email: \texttt{ksyow@upm.edu.my})
              %\email{kaisiong.yow@ntu.edu.sg}%\\
              %\emph{Present address:} of F. Author  %  if needed
           \and
            Reynold Cheng and Xiaolin Han \at
            Department of Computer Science, University of Hong Kong, Hong Kong, China
           \and
            Chenhao Ma \at 
            School of Data Science, The Chinese University of Hong Kong, Shenzhen, China
}

\date{Received: date / Accepted: date}

\maketitle

\begin{abstract}
A subgraph is constructed by using a subset of vertices and edges of a given graph. There exist many graph properties that are hereditary for subgraphs. Hence, researchers from different communities have paid a great deal of attention in studying numerous subgraph problems, on top of the ordinary graph problems. Many algorithms are proposed in studying subgraph problems, where one common approach is by extracting the patterns and structures of a given graph. Due to the complex structures of certain types of graphs and to improve overall performances of the existing frameworks, machine learning techniques have recently been employed in dealing with various subgraph problems. In this article, we present a comprehensive review on five well known subgraph problems that have been tackled by using machine learning methods. They are subgraph isomorphism (both counting and matching), maximum common subgraph, community detection and community search problems. We provide an outline of each proposed method, and examine its designs and performances. We also explore non-learning-based algorithms for each problem and a brief discussion is given. We then suggest some promising research directions in this area, hoping that relevant subgraph problems can be tackled by using a similar strategy. Since there is a huge growth in employing machine learning techniques in recent years, we believe that this survey will serve as a good reference point to relevant research communities.
\keywords{machine learning \and subgraph \and subgraph isomorphism \and subgraph matching \and community search \and community detection \and learning-based approach}
\end{abstract}

\section{Introduction}

Machine learning (ML)~\cite{Mohri_MIT2018,Alpaydin_MIT2020} is a field under artificial intelligence, which makes predictions or decisions based on past information using computational methods. It is inherently related to fields such as computer science and statistics given that data-driven learning techniques rely heavily on the availability and quality of datasets. ML is often used especially in solving some complicated multidimensional tasks that can hardly be solved by using numerical reasoning. There are many applications~\cite{Silver_NATURE2017,Mohri_MIT2018,Pugliese_DSM2021} in ML where its algorithms can be applied on various fields such as cybersecurity, financial modelling, healthcare, data governance, science and even games. As highlighted in~\cite{Pugliese_DSM2021}, ML techniques have attracted a great deal of attentions of researchers from both industry and academia in recent years. The number of publications related to ML has increased more than four times from 2016 to 2020, which demonstrates the popularity of this research direction.

Learning strategies for ML algorithms can be classified into four categories in general, which are supervised, unsupervised, semi-supervised and reinforcement learnings~\cite{Mohri_MIT2018}. The main difference between the first three learnings is on the amount of labelled data. In a \emph{supervised} learning algorithm, a machine is supplied with a set of labelled samples where the machine will learn a general pattern in mapping inputs to the desired outputs. In an \emph{unsupervised} learning algorithm, no labelled data will be supplied and a machine is required to identify structures or hidden patterns among the input data on its own. For a \emph{semi-supervised} learning algorithm, it falls in between both supervised and unsupervised learnings in which partial labelled data will be supplied in order to predict unseen points. \emph{Reinforcement learning} (RL) on the other hand is an environment-driven approach where for each state in the learning algorithm, an agent will predict over a course of actions under an environment, and learn in order to maximise some reward.

Machine learning techniques have recently been employed in addressing many graph related problems. Recall that graphs is a useful tool that is widely used in various areas such as data modelling and visualisation, in representing pairwise relationships between objects. They appear in many disciplines including mathematics, computer science, biology and sociology~\cite{Harary_ISR1953,Bondy_MACMILLAN1976,Hayes_AS2000,Pavlopoulos_BDM2011}. Graphs also involve in many real-world systems in which they can be used to model a bunch of applications in daily life such as financial markets, social networks, community networks, road networks~\cite{luo2014distributed} and information systems~\cite{DBLP:conf/edbt/LuoZXYLK23}. Many of these sophisticated systems have been addressed by using ML techniques via appropriate graph representations, demonstrating a vital role of ML techniques in graph problems.

In a mathematical context, a graph $ G=(V,E) $ is a non-linear data structure where $ V $ represents the set of \emph{vertices} (or \emph{nodes}) and $ E $ represents the set of \emph{edges} (or \emph{links}) of $ G $. A graph $ S $ is a \emph{subgraph} of $ G $, denoted by $ S \leq G $, if $ V(S) \subseteq V(G) $ and $ E(S) \subseteq E(G) $. Many interesting graph properties are generally known to be hereditary for subgraphs and minors (a substructure of a graph that is formed based on deletion and contraction operations), and many problems related to these classes of graphs have been investigated over the years.

Subgraph problems including subgraph isomorphism, subgraph counting, common subgraph, densest subgraph, community detection, community search and others are often solved by extracting structures and patterns from original graphs. Due to the complexity~\cite{Karp_1972,Garey_FREEMAN1979} of these problems and the complicated structure of certain graphs, ML methods have been developed aiming to improve the performances of traditional algorithms. Note that there are quite a number of surveys~\cite{Bengio_EJOR2021,Cappart_arxiv2021,Mazyavkina_COR2021,Peng_DSE2021,Karimi-Mamaghan_EJOR2022} focusing on ML methods in solving combinatorial optimisation and graph problems. To give an overview on learning methods that have been designed particularly for subgraph problems, in this article, we give a thorough survey by focusing on five popular subgraph problems. To the best of our knowledge, this is the first survey that reviews subgraph problems that have been addressed by learning models, where most of these models are proposed over the past three years.

We cover 21 different learning models in Section~\ref{sec:subgraph_problems}, for subgraph isomorphism counting and matching, maximum common subgraph, community detection and community search problems. We define these problems formally, and analyse the designs and performances of relevant learning frameworks. We then highlight some challenges and potential future work in the respective section. See Table~\ref{Table:Graph_problems} for more details on each of these problems.
%\liao{(1) More explanation on Table 1; (2) what is the difference between precision and effectiveness?}
%\siqiang{Maybe we can add a short discussion on the difference of our survey and previous ones.}
% More details have been included in the previous paragraph
\begin{table*}[!ht]
	\caption{A summary of subgraph problems and the respective ML solutions}
	%{\footnotesize
		%\centering
		\begin{tabularx}{\textwidth}{lp{2.5cm}p{2.3cm}p{2.3cm}ll}
			\hline\noalign{\smallskip}
			Problem & Task & Training & Algorithm & Advancement & Method (Code)\\
			\noalign{\smallskip}\hline\noalign{\smallskip}
			
			\multirow{11}{*}{\textbf{\thead{Subgraph\\ Isomorphism\\ Counting}}}
			& \multirow{7}{*}{\thead{General\\ subgraph\\ isomorphism\\ counting}} 
			& Semi-supervised & GNN & \thead{Efficiency \\ Scalability} & \thead{DIAMNet~\cite{Liu_SIGKDD2020} (\href{https://github.com/HKUST-KnowComp/NeuralSubgraphCounting}{Link})}\\
			\cline{3-6}
			& & Semi-supervised & GNN+Active Learning & \thead{Efficiency \\ Scalability} & \thead{ALSS~\cite{Zhao2021} (\href{ https://github.com/Kangfei/LSS}{Link})} \\
			\cline{3-6}
			& & Semi-supervised & GNN+Adversarial Learning & \thead{Efficiency \\ Scalability} & \thead{NeurSC~\cite{Wang2022}} \\
			\cline{2-6}
			& \multirow{4}{*}{\thead{Motif\\ isomorphism\\ counting}} 
			& Supervised & GNN & Effectiveness & \thead{LRP~\cite{chen2020b} (\href{ https://github.com/leichen2018/GNN-Substructure-Counting}{Link})}\\
			\cline{3-6}
			& & Supervised & GNN & Effectiveness & \thead{RNP-GNN~\cite{tahmasebi2021}} \\
			\cline{3-6}
			& & Supervised & GNN & Effectiveness & \thead{DMPNN~\cite{liu2022} (\href{ https://github.com/HKUST-KnowComp/DualMessagePassing}{Link})} \\
			\cline{1-6}

            \multirow{7}{*}{\textbf{\thead{Subgraph\\ Matching}}}
			& \multirow{7}{*}{\thead{Subgraph\\ matching}}
		    & Semi-supervised & Active Learning & Efficiency & \thead{ActiveMatch~\cite{Ge_BIGDATA2021}}\\
			\cline{3-6}
			& & Reinforcement & GNN+RL & \thead{Effectiveness\\ Efficiency} & \thead{RL-QVO~\cite{wang2022a}}\\
			\cline{3-6}
			& & Supervised & GNN & \thead{Effectiveness\\ Efficiency} & \thead{NeuroMatch~\cite{rex2020} (\href{ https://github.com/snap-stanford/neural-subgraph-learning-GNN}{Link})} \\
			\cline{3-6}
			& & Supervised & GNN & Effectiveness & \thead{DMPNN~\cite{liu2022} (\href{ https://github.com/HKUST-KnowComp/DualMessagePassing}{Link})} \\
% 			\cline{2-6}
% 			& \multirow{6}{*}{\thead{Matching Top\\ Subgraph}}
% 			& Supervised & GNN & \thead{Effectiveness\\ Efficiency} & \thead{IsoNet~\cite{roy2022} (\href{ https://github.com/Indradyumna/ISONET}{Link})} \\
% 			\cline{3-6}
% 			& & Supervised & GNN & \thead{Effectiveness\\ Efficiency} & \thead{Sub-GMN~\cite{lan2022a} (\href{https://github.com/zixun-lan/Sub-GMN}{Link})} \\
% 			\cline{3-6}
% 			& & Supervised & GNN & \thead{Effectiveness\\ Efficiency} & \thead{AED-Net~\cite{lan2023} (\href{ https://github.com/zixun-lan/AEDNet-Adaptive-Edge-Deleting-Network-For-Subgraph-Matching}{Link})} \\
			\cline{1-6}

			\multirow{6}{*}{\textbf{\thead{Maximum\\ Common\\ Subgraph}}}
			&
			\multirow{6}{*}{\thead{Max. common\\ subgraph}}
			& Reinforcement & Search+RL & \thead{Efficiency} & McSplit+RL \cite{liu2020} \\
			\cline{3-6}
			& & \thead{Semi-supervised,\\ Reinforcement} & GNN, Search+RL & \thead{Effectiveness} & GLSearch~\cite{bai2020} \\
			\cline{3-6}
			& & Supervised & GNN & \thead{Effectiveness\\ Efficiency} & NeuralMCS~\cite{bai2019} (\href{https://github.com/openpublicforpapers/NeuralMCS}{Link}) \\
			\cline{1-6}
			
			\multirow{9}{*}{\textbf{\thead{Community\\ Detection}}}
			&
			\multirow{9}{*}{\thead{Community\\ detection}}
			& Semi-supervised & Matrix decomposition & Effectiveness & \thead{SMACD~\cite{Gujral_SDM2018} (\href{http://www.cs.ucr.edu/~egujr001/ucr/madlab/src/SHOCD.zip}{Link})}\\ % heuristic, semi-supervised, matrix decomposition
			\cline{3-6}
			& & Unsupervised & GNN & \thead{Effectiveness\\ Scalability} & \thead{NOCD~\cite{Shchur_DLG2019}}\\
			\cline{3-6}
			& & Unsupervised & PCA & \thead{Effectiveness\\ Precision} & \thead{K-means~\cite{Chaudhary_SNAM2021}}\\ % heuristic, unsupervised, principal component analysis
			\cline{3-6}
			& & Unsupervised & GNN & Effectiveness & \thead{CE-MOEA~\cite{Sun_ITC2022}}\\
			\cline{1-6}
			
			\multirow{9}{*}{\textbf{\thead{Community\\ Search}}}
			& Interative CS & Supervised & GNN & \thead{Effectiveness\\ Efficiency} & \thead{ICS-GNN~\cite{Gao_PVLDB2021}}\\ % approximate, supervised, deep
			\cline{2-6}
			& Attributed CS & Supervised & GNN & \thead{Accuracy\\ Effectiveness\\ Efficiency} & \thead{QD-GCN~\cite{Jiang_PVLDB2022}}\\ % approximate, supervised, deep
			\cline{2-6}
			& \multirow{5}{*}{CS}
			& Semi-supervised & GNN & \thead{Effectiveness\\ Efficiency \\ Scalability} & \thead{COCLEP~\cite{Li_ICDE2023}}\\
			\cline{3-6}
			& & Supervised (Meta) & GNN & \thead{Effectiveness\\ Efficiency} & CGNP~\cite{Fang_arxiv2022}\\
			\hline			
		\end{tabularx}
		\label{Table:Graph_problems}
	%}
\end{table*}

We now list our contributions, as below:
\begin{itemize}
	\item We conduct a comprehensive review on five popular subgraph problems where learning models have been employed to solve these problems. We give an overview of these models, and discuss their designs as well as performances.
	\item We explore existing algorithms and applications for each subgraph problem. A brief discussion on their performances is also provided.
	\item We give a brief overview on learning frameworks that have been proposed to tackle some classic yet popular graph problems, where we focus on combinatorial optimisation and NP-complete problems.
	\item We list down some prospective problems that could be investigated, to either enhance further the existing frameworks or to extend them so that they can be used in solving related graph problems.
	\item This survey can be used as one of the main reference points for researchers in relevant fields (such as computer science and mathematics), particularly those who intend to employ learning-based techniques on subgraph problems.
\end{itemize}

We follow~\cite{Diestel_2010} for most of the graph terminology defined in this article, unless stated otherwise.

The outline of the rest of this article is now given. We discuss five subgraph problems in Section~\ref{sec:subgraph_problems} where we focus on subgraph isomorphism problems in Sections~\ref{sec:subgraph_isomorphism}--\ref{sec:SImatch} (we introduce subgraph isomorphism problems in Section~\ref{sec:subgraph_isomorphism}, and discuss subgraph isomorphism counting and matching problems in Sections~\ref{sec:SIcount} and \ref{sec:SImatch}, respectively). We then discuss maximum common subgraph problems in Section~\ref{sec:max_common_subgraph}, community detection problems in Section~\ref{sec:community_detection} and lastly community search problems in Section~\ref{sec:community_search}. We also give a brief overview on some classic yet well known graph problems that have been addressed by learning-based models in Section~\ref{sec:classic_graph_problems}. We concentrate on combinatorial optimisation and NP-complete problems. We then conclude this article and suggest some potential future work in Section~\ref{sec:conclusion}.

\section{Subgraph Problems}\label{sec:subgraph_problems}

We discuss five well known subgraph problems that have been addressed by ML solutions in this section. Some of them involve different tasks and they all have been tackled by different learning models.

We first introduce and define the subgraph isomorphism problem in Section~\ref{sec:subgraph_isomorphism}, before we proceed to subgraph isomorphism counting and subgraph matching problems in the two subsequent sections.

\subsection{Subgraph Isomorphism}\label{sec:subgraph_isomorphism} 

\begin{figure}[t]
	\centering
	\includegraphics[scale=1]{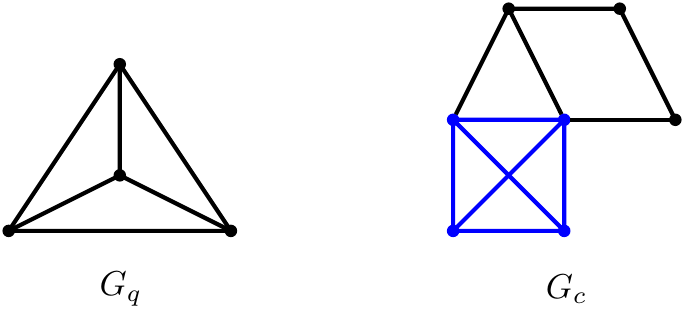}
	\caption{The graph $ G_{c} $ contains a subgraph (a complete graph of size four coloured in blue) that is isomorphic to $ G_{q} $}
	\label{fig:subgraph_iso}
\end{figure}

The \emph{subgraph isomorphism} (SI) problem takes a \emph{corpus graph} $G_c$ and a \emph{query graph} $G_q$ as input, and determines if $G_c$ contains a subgraph that is isomorphic to $G_q$. One such example is illustrated in Figure~\ref{fig:subgraph_iso}. 
	%\liao{Need to change notation in figure. Better have an example on heterogeneous graphs like~\cite{Liu_SIGKDD2020}.}

The concept of SI appears in various areas, such as in social networks (modelling and mining)~\cite{Kuramochi_ICDM2004}, databases (pattern recognition), bioinformatics (discovering protein structures)~\cite{Milo_SCIENCE2002,Alon_BIO2008}, chemoinformatics (drug development)~\cite{Huan_ICDM2003} and recommender systems~\cite{Jiang_IJCAI2017,Zhao_SIGKDD2017,davitkova2021}. 
% Many algorithms~\cite{Ullmann_JACM1976,Cordella_TPAMI2004,Carletti_TPAMI2017} based on backtracking and graph indexing have been designed in tackling the subgraph isomorphism problems. However, the search space and time grow exponentially with the sizes of the graphs due to the NP-complete nature of the problem. 

In this survey we introduce two commonly studied problems related to SI that have wide applications in practice, i.e., subgraph isomorphism counting and subgraph matching, where they both involve a high computational cost. The former gives a global estimation on the number of a query graph, while the latter determines exact vertices and edges that belong to a given subgraph. 
We discuss learning-based methods that were recently proposed in solving the two problems.
%In Section~\ref{sec:SIcount} we review the counting problem, while the matching problem will be addressed in Section~\ref{sec:SImatch}. 
% \liao{Do we need an explanation on why the SI decision problem is not included in this paper?}
\iffalse
SI decision paper:
* GIN~\cite{morris2019} and~\cite{xu2019}
*~\cite{rex2020}

minimum candidate set problem~\cite{Ge_BIGDATA2021} / SI top match:
% \subsubsection*{Learning to Discover Match by Similarity}
% ===== Interpretable Neural Subgraph Matching for Graph Retrieval
\cite{roy2022}

% ===== Sub-GMN: The Neural Subgraph Matching Network Model
\cite{lan2022a} inspired by GMN~\cite{li19d}

% ===== AEDNet: Adaptive Edge-Deleting Network For Subgraph Matching
\cite{lan2023}
\fi

% ==========
\subsection{Subgraph Isomorphism Counting}\label{sec:SIcount}

A \emph{heterogeneous graph} $ G = (V, E, L, C) $ is a labelled graph where each vertex $ v \in V $ has a label $ L(v) $ in the set of vertex label $ L $ and each edge $ uv \in E $ has a label $ C(uv) $ in the set of edge label $ C $.

%We first introduce the concept of a general \textit{heterogeneous graph} $G = (V, E, L, C)$, where each vertex $v \in V$ has a label $L(v)$ belonging to some arbitrary set of labels, and each edge $(u, v) \in E$ has a label $C((u, v))$ in an edge label set.

Given two attributed heterogeneous graphs $G_1 = (V_1, E_1, L_1, C_1)$ and $G_2 = (V_2, E_2, L_2, C_2)$, they are said to be \emph{isomorphic} if there exists a bijection $\phi: V_1 \rightarrow V_2$ such that: (1) $uv \in E_1$ iff $\phi(u)\phi(v) \in E_2$; (2) $\forall v \in V_1$, vertex label $L_1(v) = L_2(\phi(v))$; (3) $\forall uv \in E_1$, edge label $C_1(uv) = C_2(\phi(u)\phi(v))$.

We now define the \emph{subgraph isomorphism counting} (SIC) problem.
%takes a corpus graph and a query graph as input, and answers the frequency of subgraphs isomorphic to the query in the corpus graph. We consider the SIC problem for general graphs with both vertex and edge labels:
\begin{problem}
	Given a query graph $ G_q = (V_q, E_q, L_q, C_q) $ and a corpus graph $ G_c = (V_c, E_c, L_c, C_c) $, return the number of subgraphs in $ G_c $ such that those subgraphs are isomorphic to $ G_q $. 
\end{problem}
\noindent
Note that one of the variants of the SIC problem counts the number of induced subgraphs.
%For a graph $G = (V, E)$,  a \textit{vertex-induced subgraph} $G' = (V', E')$ is defined as the subgraph preserving all edges among the vertex set $V'$ with respect to graph $G$, that $\forall u, v \in V'$, edge $(u, v) \in E'$ if and only if $(u, v) \in E$. 

As determining the existence of isomorphic subgraphs in a graph is already NP-complete~\cite{Cook_STOC1971}, counting the number of subgraphs is even harder due to the exponential nature of graphs. The challenges of the problem increases as the size of a given query graph becomes larger and more complex. As a consequence, many existing solutions are only applicable to small query graphs or those with specific query patterns. 

To acquire the exact count of subgraph isomorphism, traditional methods such as enumeration and analytical approaches~\cite{Ribeiro_2022} have been employed. Enumeration methods search through the whole graph and apply categorisation algorithms, hence are not scalable when the graph or query size becomes larger. Analytical algorithms usually perform transformations on full corpus graphs and are designed for some specific patterns, which makes them difficult to apply to general subgraph queries. 
On top of the exact solutions, many approximate schemes have been developed, mainly aim to simplify and accelerate the enumeration and analytical designs. A complete introduction and evaluation on these works can be found in~\cite{Seshadhri_2019,Park_2020,Ribeiro_2022}.

% ========== subgraph counting - practice
\subsubsection*{Empirical GNNs in Subgraph Counting}

A set of researches intend to fit current deep learning models, such as neural networks, with carefully crafted structures to the SIC task for better performances, especially on large-scale and arbitrary graphs that are difficult to solve by applying traditional methods. Liu et al.~\cite{Liu_SIGKDD2020} first proposed studying the NP-complete counting problem as a graph learning task. They observed that counting queries in a corpus graph by retrieving its information is similar to the question-answering framework which is well studied in the deep learning area. 
To transfer the problem into a form that ML techniques can be applied, there are mainly two steps: (1) a representation model is applied to retrieve essential information from a corpus graph and/or a query graph; (2) a regression model then learns and predicts the count based on the input representation and queries. 
Thus, the graph information can be learned without exhaustively enumerating patterns in the full graph and the exponential term in computational complexity can be greatly reduced. By applying efficient network architectures for prediction, answering queries can also be greatly optimised. 

Representing graphs and employing neural networks to retrieve information is not a trivial task. Generally, normal neural networks learn data information that are organised in tensors. However, graph data are represented in a non-Euclidean manner, which cannot be fed directly into neural networks. 
One option is to describe an edge $uv$ as a 5-tuple $(u, v, L(u), L(v), C(uv))$ and consider it as a basic element. A graph can then be encoded as a sequence of edges, which can be learned by a series of neural networks called sequence models. 
Alternatively, a graph can be viewed as a vertex feature matrix together with an adjacency matrix. For heterogeneous graphs with edge features or for multigraphs, information can also be integrated in these matrices. A type of neural network named graph neural networks (GNNs) are proposed to specifically learn such kind of graph information. We refer interested readers to surveys such as~\cite{zhou2020,zhang2020b,wu2021} for more information. 

For the SIC model introduced in~\cite{Liu_SIGKDD2020}, several neural network structures are evaluated in respect of the performance of representing a graph. For sequence models, there are convolutional neural networks, transformer-XL and long short-term memory. Relational graph convolutional networks and graph isomorphism networks (GINs) are employed as GNN models.

The next step is to design structures that learn to predict a count for the entire corpus graph based on the encoding of corpus and query graph elements. Liu et al.~\cite{Liu_SIGKDD2020} proposed dynamic intermedium attention memory network (DIAMNet) as the interaction layer which aggregates information and outputs the prediction. Since the SIC problem involves two graphs in general, i.e., a corpus graph and a query graph, learning the graph interactions can be computationally prohibited if all combinations of vertex representations between the graphs are considered. 
The authors therefore introduced a recurrent neural network that comes with an attention mechanism, which only extracts and learns important information, in building DIAMNet. It utilises an external dynamic memory with linear complexity read-write operations to store and update query information as an intermedium. 

Liu et al.~\cite{Liu_SIGKDD2020} mainly considered a search and backtracking algorithm namely VF2~\cite{Cordella_TPAMI2004} as the traditional SIC baseline. The best and worst searching complexities of VF2 are $O(|V_c|^2)$ and $O(|V_c|!|V_c|)$, respectively, where $|V_c|$ represents the number of vertices in the corpus graph $ c $. 
In comparison, by using techniques of graph encoding and dynamic memory, the overall complexity of DIAMNet is reduced to approximately linear $|V_c|+|V_q|$, enabling its application to large-scale graphs. 

In experimental evaluation, the authors compared DIAMNet with VF2 as well as other bare sequence models and GNNs. Since existing datasets that are used in the SIC problem have limited graph and query sizes, synthetic data of different scales are generated and used alongside the realistic datasets. 
As graph convolutional operations are designed to extract the graph topological structure effectively, experiments on representation models show that graph models achieve superior performance compared to majority of the sequence models.
The model DIAMNet also outperforms pooling and attention mechanism for the interaction layer, even with moderate performances of the representation layer. 
The model empirically achieves 10--1,000 times speed-up than VF2 with acceptable errors, which is in line with complexity analysis and demonstrates the power of neural network learning in the SIC task.

% ===== A Learned Sketched for Subgraph Counting
Different from~\cite{Liu_SIGKDD2020}, Zhao et al.~\cite{Zhao2021} proposed an active learned sketch for subgraph counting (abbreviated as ALSS) model that chooses to adopt different structures for learning corpus and query graphs. The ALSS model also employs the question-answering framework, and utilises sketch learning and active learning which integrate approximate query processing in the context of relational database management system. 

Sketch learning describes the architectural design where a sketch is generated by a GNN regression model as a corpus graph representation. Inspired by the decomposition technique in traditional analytic SIC methods, query graphs are decomposed into small and basic subgraphs. Each basic subgraph corresponds to a vector of a graph representation. Hence, as the representation stage of counting a given query in a corpus graph, the learned sketch---a GIN model in ALSS---extracts the corpus graph to a series of vectorised features, each corresponds to a basic substructure in the query. 
For the module of task-specific prediction, a multilayer perceptron (MLP) with a self-attention mechanism is implemented to aggregate the corpus graph representation and the query substructures, and estimate the final count of queries. 

Typical deep learning requires a large amount of input samples, i.e., corpus graphs and query graphs with the exact count ground truth, to train a neural network. However, it is challenging in the scope of the SIC problem due to the fact that existing real-world datasets are limited in quantity and query size, and conventional exact SIC algorithms are too inefficient to generate ground truth for large-scale tasks. 
To address the problem, the authors utilised a semi-supervised learning scheme powered by active learning. \emph{Active learning}~\cite{Settles_2009} refers to a statistical ML strategy that, the algorithm can interactively perform queries to select unlabelled data and label them with desired outputs, in order to enrich the dataset for model learning and achieve comparable performances to standard supervised learning with fewer labelled samples. The active learner, which is an uncertainty sampling algorithm, queries the sketch model itself to choose the query graph for the following learning scenarios. 

For experiments, the authors adopted recent benchmark datasets and algorithms from~\cite{Park_2020}. The datasets are real-world data on social network and web graphs with query size up to 32. Two baselines are WJ~\cite{Li_SIGMOD2016} and IMPR~\cite{Chen_tkdd2018} that perform SIC estimations based on random walk samplings. 
Experiments show that traditional algorithms suffer from severe underestimation due to sampling failure for queries with at most eight vertices, and are not applicable to larger queries. In comparison, ALSS performs counting that is close to ground truth consistently for both small and large queries. 
ALSS also surpasses WJ and IMPR in efficiency, with 1--2 order of magnitude faster in query time especially for large corpus and query graphs. 
% The ALSS~\cite{Zhao2021} model combines sketch learning and active learning, which integrates approximate query processing (AQP) in the context of relational database management system (RDBMS). The ALSS model is also under the question-answering framework, where a GIN is used to generate a learned sketch of the corpus graph by decomposing it into smaller substructures for training, and an MLP with self-attention is the prediction module that aggregates representations and estimates the count of queries. To improve the generalisation of the supervised sketch learning, the authors further employed uncertainty sampling as active learning while learning the graph sketch. 

% ===== Neural Subgraph Counting with Wasserstein Estimator
Wang et al.~\cite{Wang2022} further addressed the efficiency and scalability issues of the previous ML solutions~\cite{Zhao2021} for the SIC problems, while providing more effective designs in retrieving graph information and conducting semi-supervised learning. They proposed neural subgraph counting (NeurSC) that consists of an extraction module to extract representation differently from corpus and query graphs, and an estimator that integrates representation and output count predictions. 
They observed that not all graph structures are useful with regard to counting subgraphs. Graph representation structures of existing ML approaches such as DIAMNet and ALSS are inefficient since they need to process redundant information. The extraction module of NeurSC is hence designed to extract relevant vertices from a corpus graph based on vertices in given query graphs. In other words, it prunes the large graph to filter out candidate substructures instead of presenting the entire corpus graph as a whole, thus a representation can be generated on both corpus graph substructures and query graphs without prohibitive time overhead. The scheme is also helpful for avoiding the influence from unpromising vertices in the graph, thus produces a simple but effective representation. 

The estimator then employs dual GNNs to generate subgraph representations and learn vectorised features. Specifically, an intra-graph network, which is implemented as a GIN, captures the information individually for query graphs and corpus graph substructures. Another GNN, the inter-graph graph attention network (GAT)~\cite{Velickovic_ICLR2018}, is exploited to specially process the interrelationship between query and corpus graphs. This is achieved by constructing the correspondence of vertices in query graphs and those in corpus graph substructures and a bipartite graph, and the inter-graph GNN learns features from the bipartite graph. The features acquired from the two GNNs are then aggregated and passed to an MLP for predicting the final count. 

The authors noticed that, it is not guaranteed that the representation generated by GNNs between the query graphs and corpus graph structures are similar, even for the same patterns. To cope this, a Wasserstein discriminator is proposed to conduct adversarial training. The discriminator minimises the Wasserstein distance of corresponding query and corpus substructures in the representation space, so as to maintain the exactness of the extraction algorithm and the dual GNNs. 

To evaluate the power of NeurSC, they compared it with conventional SIC algorithms~\cite{Park_2020} as well as ML methods DIAMNet~\cite{Liu_SIGKDD2020} and ALSS~\cite{Zhao2021} on seven realistic datasets. Conventional benchmarks present similarly as in~\cite{Zhao2021} that sampling-based algorithms meet sampling failure especially for large and complex graphs, and summary-based ones suffer from underestimation due to their prior assumptions. 
In general, learning models outperform traditional models, and NeurSC consistently outputs better predictions than ALSS, while DIAMNet exceeds execution time limit on large datasets. 
The authors also particularly studied the effectiveness of NeurSC and ALSS on graphs with various complexity, density, and diameter. Both methods tend to perform better on graphs that have lower complexity, larger density and smaller diameter. The NeurSC model consistently achieves better accuracy, and its advancement is observed to be greater for large and complicated graphs. 
When considering efficiency, it is the fastest for the overall query time, though ALSS is faster in training as it processes less graph information, while NeurSC requires the expensive adversarial training.

% ========== subgraph counting - theory
\subsubsection*{Expressive GNNs in Motif Counting}

In solving the SIC problem, there exist another line of researches that are similar to the traditional analytical approaches, where specific motifs are examined and theoretical guarantees of solutions are provided. These works are based on recent advances of the expressiveness of GNNs, and models learn the task themselves instead of adding carefully designed components. 
In a more general context of studying the expressive power of learning-based algorithms, it is proven~\cite{morris2019,xu2019} that GNNs based on neighbourhood-aggregation schemes can be as powerful as the Weisfeiler-Lehman (WL) test~\cite{Weisfeiler1968,ARVIND202042}, which is effective in most cases to solve the SI problem. Such a discovery encourages researchers to apply GNNs in graph problems, especially on variants of the SI problem, pursuing the full power of learning algorithms in effectiveness and efficiency. 

% ===== Can Graph Neural Networks Count Substructures?
By extending this framework, Chen et al.~\cite{chen2020b} studied the expressiveness of GNNs specifically in the SIC problem, which mainly includes the two common structures: (1) a message passing neural network (MPNN) that computes graph representations by iterative transformation and local aggregation of vertex and edge features, including a wide range of popular models namely graph convolutional network (GCN)~\cite{Kipf_ICLR2017}, GIN~\cite{xu2019} and GraphSAGE~\cite{Hamilton_NIPS2017}; (2) an invariant graph network (IGN) that represents graph information by high-order tensors and processes them by permutation equivariant functions. It is known that a 2-order IGN (2-IGN) can approximate an MPNN well in expressiveness. 

Chen et al.~\cite{chen2020b} proved that in regard to counting induced query graphs with at least three vertices, the MPNN and 2-IGN cannot count all subgraphs accurately as they fail to distinguish certain graphs. However, these models are powerful enough to count subgraphs for star-shaped patterns. Intuitively, this is because the local aggregation scheme of these architectures gradually grows a computation tree when calculating representations. However, essential structural information such as vertex identity and connectivity are lost. To overcome this expressiveness bound, more powerful schemes such as a high-order $k$-IGN test need to be applied. 
They proved that $T$ iterations of $k$-WL (equivalent to $k$-IGN) can be expressive on queries of scales up to $(k+1) \cdot 2^T$ vertices, providing a perspective of explanation of the ability of deep and iterative GNNs. 

Based on the expressiveness analysis, Chen et al.~\cite{chen2020b} proposed a type of GNN architecture, namely local relation pooling (LRP), that is empirically evaluated to be powerful in solving the SIC task. The concept is based on the observation that substructures present themselves in local neighbourhoods. As both MPNN and 2-IGN are bounded by their iterative equivariant aggregation, the authors bypass this with relation pooling to effectively learn the subgraph information. 
The LRP is designed to improve the 2-IGN aggregation that, instead of applying global permutations on the whole graph, not necessarily permutation-invariant functions are computed on the $l$-hop neighbourhood of each vertex in the graph. Then the neighbourhood information are aggregated by pooling to form the global representation of the entire graph. Such mechanism enables efficient and less constrained computation on the local permutation, while being specially optimised to retrieve patterns of a radius less than $l$. When stacking LRP with more than one layer, Deep LRP is obtained and it quickly scales up in the ability to process large queries. 

The authors evaluated Deep LRP by using both synthetic random data and realistic molecular data, comparing to baselines of typical MPNN structures and 2-IGN. The queries are motifs with a maximum order of four. Experimental results show that for the MPNN and 2-IGN, their counting accuracy is better for counting star-shaped patterns than other shapes (such as triangles), which is in line with the theoretical analysis on MPNN expressive power. 
For LRP and Deep LRP, the performance is among the best regardless of the shape of motifs, indicating that the LRP structure has the expressive power to count motifs as expected. 

% ===== Counting Substructures with Higher-Order Graph Neural Networks: Possibility and Impossibility Results
Given that the SIC problem is used as a criterion of evaluating the power of GNNs~\cite{chen2020b}, more studies in providing relevant analyses and designing GNN schemes are found. A recent work gives a general bound regarding the GNN expressive power in performing the SIC task~\cite{tahmasebi2021}. 
% To represent an arbitrary corpus graph $G_c$ with $n$ vertices, representation functions like GNNs use an $s$-value scheme to encode $t$ subgraphs of $G_c$. Then, to perform correct counting with query graph with $k$ vertices, the number of subgraphs need to be encoded is $t = \tilde\Omega (n^{\frac{k}{s-1}})$, where $\tilde\Omega (\cdot)$ is $\Omega (\cdot)$ up to poly-logarithmic factors. 
It shows that to perform correct counting with a corpus graph of $n$ vertices and a query graph of $k$ vertices, the complexity of any representation functions like GNNs is $n^{\Omega(k)}$ in the worst-case scenario. More strictly, even for weakly sparse graphs with the average degree $O(1)$, no representation algorithm can complete the count in linear time. 

In the paper, Tahmasebi et al.~\cite{tahmasebi2021} developed the expressiveness theory especially on the scheme of a recursive neighbourhood pooling GNN (RNP-GNN), which is more flexible than IGN and LRP that are computationally expensive at high orders. Compared to LRP, instead of calculating the $l$-hop neighbourhoods of a vertex $v$, it represents vertex-wise local information based on a collection of subgraphs $\mathcal{N}_l (v)$ with radius smaller than $l$. 
Such a representation can be computed in a recursive manner that, to calculate the local representation of vertex $v$ in a graph $G$, it can be aggregated by the representation of the set of neighbouring vertices $G_v = \mathcal{N}_l (v) \setminus \{v\}$. The recursion is beneficial in reducing the size and hierarchy of subgraphs, which simplifies the computation while providing sufficient expressiveness. 
Notably, RNP-GNN does not require other encoding schemes, as recursive pooling itself is sufficient when applied recursively for the SIC problems.

The authors proved that, by appropriate recursion of $k$ computations, RNP-GNN has the expressive power to perform counting on any induced and non-induced subgraphs with $k+1$ vertices. More generally, it is able to approximate any arbitrary graph functions related to subgraphs of size $k+1$. 
The complexity of RNP-GNN is highly related to the maximum degree of a graph. For any corpus graph of size $n$, the worst-case complexity of RNP-GNN is $n^k$, which is the same as $k$-WL, compared to LRP that is not bounded. For extreme cases of weakly and strongly sparse graphs, RNP-GNN achieves a linear time complexity, which is close to the theoretical complexity bound mentioned above. 
%For extreme cases of weakly ($\bar{d} = O(1)$) and strongly sparse graphs ($\max(d) = O(1)$), RNP-GNN can achieve $O(n)$ and linear complexity, which is close to the theoretical complexity bound above. 

% ===== Graph Convolutional Networks with Dual Message Passing for Subgraph Isomorphism Counting and Matching
Suppose $ G $ is a graph, its \emph{line graph} $\mathcal{L}(G)$ is obtained as follows: (1) every vertex of $\mathcal{L}(G)$ corresponds to an edge of $G$, and (2) two vertices in $\mathcal{L}(G)$ are adjacent if the two associated edges in $G$ share a common endvertex. In order to apply duality in the SIC problem, Liu and Song~\cite{liu2022} chose to improve the MPNN scheme by designing dual message passing neural networks (DMPNNs) in exploiting graph duality. It is known that there is a one-to-one correspondence between isomorphisms of simple graphs and isomorphisms of their dual graphs~\cite{whitney1992congruent}. They proved that such a correspondence holds for connected directed heterogeneous multigraphs with reversed edges. Hence learning to count isomorphism on raw corpus graphs is equivalent to modelling line graphs for an algorithm.
%A \textit{line graph} $\mathcal{L}(G)$ is achieved through edge-to-vertex transformation on a graph $G$, that each edge in $G$ corresponds to a vertex in $\mathcal{L}(G)$, and every two edges in $G$ that have a common vertex is represented by an edge between their corresponding vertices in $\mathcal{L}(G)$.

To fully make use of both vertex and edge information in a graph, DMPNN conducts dual message passing in generating graph representations. In the message passing stage of conventional MPNNs, the vertex state is iteratively computed based on a graph convolution operation of the previous state and graph adjacency, where the graph connectivity is only included in its adjacency. 
In DMPNN, it becomes a dual computation where the vertex and edge features are represented by two matrices, and are updated by the integration of the previous information. Both vertex and edge information are explicitly retrieved and can be learned by the model parameters. Hence, its expressive power in representing graph substructures can be improved. 
The dual message passing scheme is applicable on both heterogeneous and homogeneous graphs. 

The DMPNN model can also be integrated with LRP that acts as a pooling layer after message passing, which is named as DMPNN-LRP. The authors evaluated both DMPNN and DMPNN-LRP on small queries based on settings similar to LRP, as well as large queries with at most eight vertices using more general datasets. 
Comparisons against sequence models and MPNNs show that DMPNN outperforms them in the SIC task with respect to the root mean square error and mean absolute error, as the merit of its utilisation of both vertex and edge features. 
Note that LRP alone performs poorly for large queries. However, DMPNN-LRP surpasses both DMPNN and LRP in both small and large queries, though the running time is almost doubled. 
It is also demonstrated that by including additional reversed edges information on heterogeneous graphs, models with graph convolution perform better in a consistent manner. The improvement indicates the importance of learning complex local graph structures in solving subgraph problems.

% ===== future work
\subsubsection*{Future Direction}

% Given the effectiveness of this technique, one could consider to apply it on other NP-complete problems with appropriate modifications. Since the learning framework gives a promising result in the subgraph isomorphism counting problem, it could also be extended or generalised particularly to other variants of the SI problem.
As the power of learning-based techniques, especially the MPNN model, has been explored theoretically and empirically, the SIC problem becomes a promising topic that demonstrating the effectiveness of learning frameworks in subgraph problems. Works reviewed in our survey suggest that, SIC models using GNNs generally outperforms traditional methods regarding precision, efficiency and scalability. The expressive power of GNNs is also proved to be as powerful as the WL test. Given that SI and SIC problems are becoming widely recognised criteria of GNN expressiveness, more studies can be conducted to explore the power of GNNs in these problems from various perspectives, such as the approximation ratio that was investigated in other combinatorial problems~\cite{sato2019}. 

We notice that for current works, empirical GNN solutions usually contain handcrafted structures and are more efficient in general tasks, while GNN architectures stressing expressive power come with theoretically guaranteed efficacy and can perform learning automatically, but can only apply to specific queries and often have less competitive efficiency. 
In the future, more efforts can be put into bridging the theoretical guarantee and empirical performance, to obtain a more general solution for the SIC problem and its variants, and even other NP-complete problems. 
Additionally, the SI problem can be related to other problems through conversions, such as the maximum common subgraph problem that will be discussed in Section~\ref{sec:max_common_subgraph}. It is also possible that algorithms for these problems can be integrated.

\subsection{Subgraph Matching}\label{sec:SImatch}

The subgraph matching (SM) problem is known to be related to the SI problem on a corpus graph and a query graph, which is formulated as follows: 
\begin{problem}
	Given two attributed graphs, query $ G_q = (V_q, E_q, L_q, C_q) $ and corpus $ G_c = (V_c, E_c, L_c, C_c) $, determine all subgraphs of $ G_c $ that are isomorphic to $ G_q $.
\end{problem}

Different from the SIC problem that only returns the global count estimation of corpus-query graph pairs, which can be achieved without searching the entire graph exhaustively, the SM problem needs to find the exact graph elements, i.e., vertices and edges, that satisfy the bijection to the query graph. Despite its simple definition, it is more difficult in general as its search space greatly increase if a corpus graph has a large set of candidate vertices. In current benchmark datasets~\cite{Ge_BIGDATA2021}, the total SI count can go up to approximately $ 10^{384} $. 

Traditionally, constraint propagation~\cite{Lee_VLDB2012,Luo_SIGMOD2020} and tree search are common techniques used in solving the SM problem.  
Constraint propagation~\cite{Moorman2021} lists potential corpus vertices that could match with each query vertex, and possible matches can be reduced by applying local constraints.
Tree search algorithms~\cite{Carletti_GBRPR2015,Carletti_GBRPR2017} on the other hand enumerate vertices in a query graph and attempt to match them in a corpus graph using a search tree. It retains a search state while traversing the tree, and backtracks whenever it lands at the end of a branch. The order of matching vertices is critical as different search and backtrack sequences can vary greatly in the generated search tree sizes~\cite{Luo_SIGMOD2020}.
To solve the SM problem, these two methods can also be merged.

\subsubsection*{Learning to Optimise Search}

ML approaches for the SM problem have been addressed by Ge and Bertozzi~\cite{Ge_BIGDATA2021} using active learning for matching solution optimisation. Similar to the SIC problem, datasets in the SM problem contain insufficient labels while the original solution space is too large for exact algorithms to generate ground truth. Hence, semi-supervised learning is required to specifically select query vertices to obtain additional information for SM algorithms so that less labelling work are required. 
As a semi-supervised solution, the optimisation goal of active learning in~\cite{Ge_BIGDATA2021} is to provide a final solution to the SM tree search by eliminating extraneous isomorphic subgraphs, hence the solution space will have a modest solution count that is suitable for the SM algorithm. 

To employ active learning in solving the SM problem, the goal is to identify and choose query vertices that could potentially reduce the search space. For this purpose, three strategies have been proposed where the authors focus on the degree centrality measure of query vertices, sum of candidates for adjacent vertices and edge entropy. Intuitively, the first two strategies choose vertices that are most \emph{relevant} to other vertices, while the last simplifies the complex part of a graph as edge mapping is usually more complex and contains richer information than vertex mapping in the SM problem. These techniques are applied after potential candidates are identified by using constraint propagation, to provide real-time candidates for active learning labelling. 

They conducted three case studies based on datasets that cover domains of economic transactions, social media and transportation. For example, on the IvySys graph with $2,488$ corpus vertices, each query vertex would have typically more than $2,000$ corresponding candidates, resulting the search space to become as large as $10^{100}$. By applying active learning criteria, the number of candidates per vertex can be significantly reduced to less than $100$. It even filters out the only candidate in the corpus graph for certain vertices, which can greatly simplify the search. 
Comparisons on the three template identifications suggest that the optimal choice of criteria differs among the datasets and depends heavily on the graph topology and structure. The authors suggested that additional metrics that involve the graph topology to be investigated.

% ===== Reinforcement Learning Based Query Vertex Ordering Model for Subgraph Matching
Wang et al.~\cite{wang2022a} used a similar idea on learning scheme to optimise backtracking search matching order by RL algorithms. 
%Recall that RL describes a set of learning methods that defines an environment and an agent that interacts and learns from the environment to reach a solution.
Applying RL to combinatorial optimisation problems has been of particular interests in recent ML communities~\cite{Mazyavkina_COR2021}. 
In the context of graph problems, the \emph{state} in RL usually refers to the current partial solution to the overall problem (e.g., the subgraph of matched vertices in the SM search), and an \emph{action} is to change the solution (e.g., selecting a new query vertex for matching). 
Usually, RL methods use a \emph{reward} function to indicate how an action is carried out based on the current state changes towards the optimisati on goal. 

The authors stated two challenges in the existing solutions for the SM problem. First, matched vertices are represented and selected only based on the local neighbourhood information such as connectivity. The topology from corpus and query graphs are not fully utilised, providing possibility of overlooking certain structural and attribute information. Second, selection strategies are usually static and greed heuristics, which imply that they cannot fit to the global procedure and may be trapped in suboptimal solutions. 
Considering the above constraints, they proposed a reinforcement learning based query vertex ordering (RL-QVO) model for subgraph matching using the backtracking strategy. 

To address the first problem, the RL-QVO model employs the GCN~\cite{Kipf_ICLR2017} architecture to generate representations from both corpus and query graphs. As these graphs do not contain edge labels, the features are designed to include vertex degree, vertex label and query vertex encoding. An RL policy search with the GCN as the policy network is then designed to conduct learning for matching order. The state of RL contains two feature matrices for the sequence of selected vertices and the entire query graph, as each vertex is expressed as a feature vector generated by the GCN. For each action, the algorithm selects a new vertex to add into the sequence based on the predicted probability. 
The reward is the key part of adopting an RL algorithm. For RL-QVO, the reward function is the weighted sum of: (1) enumeration reward that pays more attention to complex queries and therefore simplifies the search; (2) validation reward that checks the probability validity of model selections; (3) entropy reward that encourages output complicated probability distribution. 

The experiment evaluation of RL-QVO is conducted against five conventional SM algorithms, covering search, filtering and hybrid approaches, on six real-world datasets that include a million-scale corpus graph~\cite{Luo_SIGMOD2020}. It is shown that RL-QVO is able to find matches with up to two orders of magnitude faster than other baselines. By the power of RL alongside, its search performance is particular superior on large corpus and query graphs, and in finding those hardest matches in the same graph. 
It is also evaluated that RL-QVO benefits from fast and effective training, and the overhead of training and inference of the ML module is neglectable compared to the enumeration search.

% ===== 
\subsubsection*{Learning to Match in Embedding Space}

Due to a significant increment of ML expressive power in recent advances, another approach to solve the SM problem is to perform the matching directly using graph representations generated by neural network models in an end-to-end fashion. We present two relevant studies under this topic. However, these works do not perform the exact same SM task as those discussed earlier, mainly due to the limited designs of their matching strategies. In other words, they do not return individual vertices and edges that matched in subgraphs, but some alternative results.
%Nonetheless, we think they are inspirational in elaborating the application of learning embedding in SM. 

% ===== Neural Subgraph Matching
Ying et al.~\cite{rex2020} first applied the idea of utilising neural network learning for matching subgraphs in the representation space by proposing a model namely NeuroMatch. Similar to SIC empirical models, it involves a two-stage learning procedure: (1) learn a representation model to process corpus graphs and query graphs, either separately or jointly; (2) implement a prediction model to find matches in the representation space. 
% Note that different from SIC regression module which only returns count estimation, the SM prediction need to contain all information of the specific match. 
Note that NeuroMatch is different from most SM search algorithms, where it requires the utilisation of the anchor technique to select a unique vertex from a given graph as the anchor and then perform matching based on it. Lan et al.~\cite{lan2023} commented that this may cause underestimation on matches in certain cases. 

The NeuroMatch representation is designed separately for corpus and query graphs. In the former, it borrows the concept of substructure decomposition to enumerate each vertex and create the $l$-hop neighbourhood as the substructure. In the latter, an anchor vertex is selected for representing the query. 
A GNN is then trained to embed the substructures and queries into the representation space by exploiting order embedding~\cite{McFee2009}, where the subgraph relationship can be measured directly. In practice, corpus graphs can be decomposed and processed beforehand to cut the computation cost. 

For matching the corpus substructures and queries in the embedding space, NeuroMatch chooses to adopt a simple voting algorithm. Given a query graph, NeuroMatch enumerates anchor vertices in the corpus-query pair and compares the representation of their neighbourhood. If by order embedding, the similarity between the anchor pair is above a certain threshold, it counts a valid vote for the subgraph match. However, a natural limitation of this scheme is that, it can only decide if a query is isomorphic to a substructure, rather than giving the exact match of vertices and edges of the isomorphism. 

Due to the limited ability of NeuroMatch, it can only be evaluated based on particular tasks to infer the isomorphism existence in corpus graphs, and measured by prediction accuracy and area under the receiver operating characteristics. For experiments, they found that the model is up to 100 times faster compared to exact SM algorithms~\cite{Cordella_TPAMI2004} and 18\% more accurate than approximate methods~\cite{lu2016} with similar runtime. 

% ===== Graph Convolutional Networks with Dual Message Passing for Subgraph Isomorphism Counting and Matching
The DMPNN model~\cite{liu2022} introduced in Section~\ref{sec:SIcount} is also applicable to the SM problem since it is proposed as a general GNN architecture for multiple purposes. Instead of predicting a global count in the SIC task based on vertex and edge representations, it outputs the number of occurrences in isomorphic subgraphs for every corpus vertex in the SM task. The authors claimed that such local results verify the model ability for subgraph matchings. 
Its learning scheme also follows the process of first embeds structural information of vertices and edges, then learns to predict vertex-wise regression on SM occurrences by higher model layers. 

For evaluating matching performances, the authors utilised the average graph edit distance, which is a common metric in graph matching related scenarios. Experimental results are similar to those in the SIC task, where both DMPNN and DMPNN-LRP are more effective than conventional sequence models and MPNN designs, especially for complex corpus-query pairs. The accuracy improvement of the structure is brought by the explicit consideration of edge information.

% ===== future work
\subsubsection*{Future Direction}

Similar to the SIC problem, the SM problem has a wide relation with other graph problems in NP, such as graph matching and subgraph similarity. The DMPNN model~\cite{liu2022} typically proved to have a general expressiveness in performing both the SIC and SM tasks, which shows certain relation for these two SI-related problems underneath. However, as we stated, current approaches of utilising GNN learning to solve the SM problem are limited by their output prediction. Hence, a possible future direction is to improve the ability of GNN structures to make exact inference on matched components. 
There are also a set of learning-based efforts on solving a similar SI candidate set problem (i.e., finding \emph{any} SI matching) based on similarity measures, which may provide another perspective of view~\cite{roy2022,lan2022a,lan2023}. 

% importance of edges~\cite{liu2022,roy2022,lan2023}

\subsection{Maximum Common Subgraph}\label{sec:max_common_subgraph}

Instead of searching query graphs in a corpus graph, the \emph{maximum common subgraph} (MCS) problem aims to determine isomorphic subgraphs in two parallel graphs. Specifically, it targets at finding the largest common subgraph based on two input graphs~\cite{Bokhari1981,BUNKE1997689}. 
\begin{problem}
	Given two input graphs $G_1$ and $G_2$, the maximum common subgraph problem finds a maximum induced subgraph $ H $ such that $ H \subseteq G_{1} $ and $ H \subseteq G_{2} $.
	%Given two input graphs $G_1 = (V_1, E_1)$ and $G_2 = (V_2, E_2)$, the maximum common subgraph problem finds a maximum induced subgraph $G_1' = (V_1', E_1')$ and $G_2' = (V_2', E_2')$ such that $G_1'$ and $G_2'$ are isomorphic. 
\end{problem}
\noindent
A common variant of the MCS problem is defined based on two labelled graphs $G_1 = (V_1, E_1, L_1)$ and $G_2 = (V_2, E_2, L_2)$. Hence, the problem requires a bijection $\phi$ between the vertex sets $V_1' \subseteq V_1$ and $ V_2' \subseteq V_2$ of two subgraphs such that $\forall v \in V_1'$, we have $L_1(v) = L_2(\phi(v))$. 
Another variant namely the maximum common connected subgraph problem on the other hand requires the induced subgraph being connected.

Having an intuitive definition, the MCS problem is widely used as a measure to represent graph similarity. It has a broad utilisation in various domains such as information retrieval~\cite{yan2005,caoyangwang2011}, programme analysis~\cite{Djoko1997,park2011} and biochemistry~\cite{faccioli2005,COOTES20071126,Giugno2013}. 

The MCS problem and many of its derived variants are NP-hard and computationally challenging. Before learning-based methods emerge, traditional efforts have been developed based on constraint programming~\cite{vismara2008,ndiaye2011,McCreesh2016} and integer programming~\cite{BAHIENSE20122523}. The state-of-the-art algorithm McSplit~\cite{mccreesh2017} utilises a branch and bound (BnB) heuristic search algorithm, which is able to exhaust the search space efficiently. For each iteration, the algorithm grows a candidate subgraph by selecting vertices at its branching point, then determines whether to keep or prune the current state. However, it requires exponential time in the worst-case scenario, which prohibits its application to large-scale graphs.

\subsubsection*{Learning to Optimise Search}

% ===== A Learning Based Branch and Bound for Maximum Common Subgraph Related Problems
To apply learning methods to the MCS problem, Liu et al.~\cite{liu2020} adopted the way of setting an RL algorithm alongside the search to achieve a better performance. They proposed McSplit+RL that enhances the search heuristics McSplit~\cite{mccreesh2017} by introducing RL strategies during the search process.

Applying RL alongside traditional models for the MCS problem is according to the consideration that a search algorithm can be regarded as an agent that performs a sequence of actions to find the MCS~\cite{liu2020}. For the BnB algorithm, each branching choice can be seen as an action. 
Hence, the RL strategy can be applied alongside the BnB algorithm to learn and choose vertices in growing a candidate subgraph. Different from the traditional McSplit algorithm that selects vertices at branching point simply based on vertex degrees, the RL algorithm can learn to choose the optimal vertex. 

Designing reward function is critical in exploiting RL algorithms, as it significantly affects the overall behaviour of the agent. The authors refer to the principle of BnB in which its goal is to reach a leaf of the search tree as early as possible to minimise the size of the search space for subsequent steps. Hence, the reward is designed to lead the algorithm in completing the search branch easily.

Performance evaluations show that more instances can be solved successfully by McSplit+RL with less search time compared to the ordinary McSplit algorithm. 
In investigating factors where McSplit+RL outperforms the baseline algorithm, the authors found that it involves more vertices during the search, which indicates that RL enables a more diverse search space and is hence useful in enhancing BnB search.

\subsubsection*{Learning to Directly Optimise Solution}
% ===== GLSearch: Maximum Common Subgraph Detection via Learning to Search
% It adopts GNN-based Deep Q-Network (DQN) as its RL model, and the goal of its RL agent is to directly maximise the common subgraph size. By directly optimising the problem and output MCS candidates, the algorithm is able to find solutions quickly and effectively, without much backtracking and pruning, and prevent trapped in local optima. 

Bai et al.~\cite{bai2020} proposed an RL-powered model, namely GLSearch that is built on top of BnB search under the learning to search framework. Their framework differs from~\cite{liu2020} by deeply reshaping the process of reaching the optimal solution. It demonstrates that ML models are able to solve the MCS problem directly by learning rather than assisting conventional search algorithms. 

The GLSearch model applies RL to aid BnB search, but primarily adopts the deep Q-network (DQN) model, which learns based on a quality function $Q$ describing the current state-action combination, as its RL backbone. It provides the model with the ability to conduct semi-supervised learning where the model is first pre-trained with sole BnB search on small graphs, then transfer and perform predictions on large datasets. 

The optimisation goal of utilising RL in the MCS problem is also altered. Different from~\cite{liu2020} that aims to simplify the BnB search, the goal of DQN application in GLSearch is to directly optimise the problem outputs, i.e., to maximise the candidate common subgraphs. This strategy is able to enhance the search process on a higher level. The RL algorithm not only learns to perform branching point choices, but also adjusts the search order to find solutions quickly and effectively without much backtracking and pruning, and prevent traps in local optima. 

To adapt the RL optimisation, the $Q$ function that integrates both the current state and action as well as guides the agent, is carefully designed to exploit the graph embedding produced by a GNN model. The $Q$ function includes vertex neighbouring information as well as graph information of the two input graphs, in order to capture both local and global information for better representations.  

Empirical evaluations show that the GLSearch model produces better results comparing to McSplit and McSplit+RL, under the same computational resources. The improvements are significant especially on large-scale graphs, where the baselines usually fail to reach solutions in early stages due to their searching heuristics.

% ===== Neural Maximum Common Subgraph Detection with Guided Subgraph Extraction
Different from~\cite{liu2020} and~\cite{bai2020}, another promising line of researches aims to apply learning-based models directly onto solutions of the MCS problem. 
Bai et al.~\cite{bai2019} firstly proposed the NeuralMCS model to apply deep learning approaches directly to the MCS problem without adopting conventional searches alongside. 
Its idea is to apply a GNN network to learn a vertex-wise similarity matrix between two input graphs, and gradually generate the candidate subgraphs based on the similarity metric. 
However, the current form of such approach is highly constrained for demanding a large amount of precomputed MCS instances as training samples, which limits its practical applications.

% ===== future work
\subsubsection*{Future Direction}

For future works, the application of learning algorithms on the MCS problem can be of particular interests for both ML and graph communities. 
First, the MCS problem is highly related to many subgraph problems, as graph isomorphism and SI problems (Section~\ref{sec:subgraph_isomorphism}) can both be considered as its special forms. Studies have also shown that the MCS problem can be converted to the maximum clique detection~\cite{levi1973,McCreesh2016} problems. Hence, existing solutions can potentially be transformed to these related problems in NP. 
Second, the pipeline of bridging learning-based and conventional approaches can be insightful when solving similar graph problems. In general, traditional heuristic search techniques are usually mature and provable, but less flexible. Fitting these algorithms to learning frameworks can bring benefits such as in improving efficacy or scalability.
% \siqiang{This subsection about Max Common Subgraph is quite short. Given that this problem is one of 4 main problems we aim to discuss, we should discuss them in much more details, even to algorithmic steps.}

\subsection{Community Detection}\label{sec:community_detection}

\begin{figure}[t]
	\centering
	\includegraphics[scale=0.8]{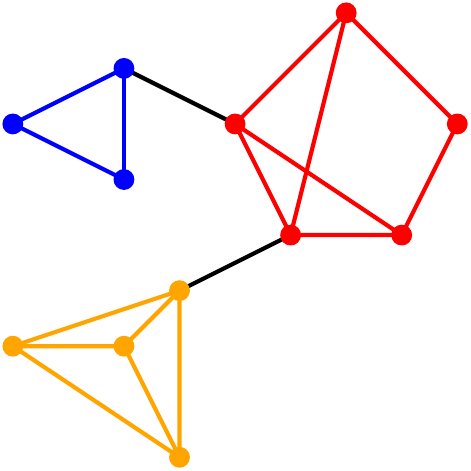}
	\caption{A community that has three clusters, coloured in three different colours}
	\label{fig:community_detection}
\end{figure}

Community detection (CD)~\cite{Fortunato_PR2010} is useful in identifying relationships among vertices and understanding the structure of a given graph. Its main objective is to partition a graph $ G $ into different groups in a way such that the set $ U $ of vertices within a group are densely connected among themselves, and at the same time sparsely connected with vertices in $ V(G) \setminus U $. One such example is shown in Figure~\ref{fig:community_detection} where a community is partitioned into three groups.
\noindent
There exist numerous applications in CD where it is beneficial in fraud detection~\cite{Pinheiro_SUGI2012} and in identifying false or missing links in a network. Elements with identical characteristics can also be detected in social networks analysis~\cite{Girvan_PNAS2002}.

Various frameworks have been developed in solving CD problems for different types of networks. Most of the traditional works concentrate on just the two basic elements, i.e., vertices and edges, of a given graph~\cite{Leskovec_WWW2010}. Some conventional methods~\cite{Fortunato_PR2010,Su_TNNLS2022} in solving the problem include graph clustering, hierarchical clustering, spectral clustering and optimisation. These methods usually detect communities according to structures of graphs. On the other hand, some algorithms merge different aspects of a graph to enhance the accuracy and preciseness of outputs in CD problems~\cite{Tang_ICDM2009,Papalexakis_ICIFS2013}. Due to the high dimensional features and complex topology of certain graphs, learning-based frameworks for CD problems have been developed. Recent findings~\cite{Shchur_DLG2019} also reveal that learning-based algorithms are more accurate and scalable. A relevant survey can be found in~\cite{Su_TNNLS2022}.

In this section, we focus on four learning frameworks in solving CD problems that are related to (1) overlapping communities, (2) attributed networks and (3) Covid-19 dataset.

\subsubsection*{Overlapping Community Detection}

Real complex networks often consist of overlapping communities~\cite{Yang_ACM2014}, and it is essential to develop models to partition such complex networks in an effective and efficient manner. Over the years, several non-learning methods have been proposed in discovering overlapping communities. We first give an overview on three of the proposed methods: probabilistic inference, heuristics and non-negative matrix factorisation (NMF). A detailed survey for overlapping CD problems can be found in~\cite{Xie_ACM2013}.

\emph{Probabilistic inference}. Latouche et al.~\cite{Latouche_AAS2011} proposed a random graph model that generalises the stochastic block model~\cite{Nowicki_JMSA2001} to detect overlapping clusters in complex networks. They focused on directed graphs, and showed that the model is identifiable within equivalence classes. Zhou~\cite{Zhou_AIS2015} designed an edge partition model in partitioning undirected relational networks by using hierarchical gamma process. The model can be used for link prediction and overlapping CD, and is scalable to networks that have thousands of vertices. Todeschini et al.~\cite{Todeschini_JRSS2020} generalised existing probabilistic models in dealing with overlapping community structures for simple graphs (loops are allowed), which can also be extended to bipartite graphs. An algorithm based on Markov chain Monte Carlo methods was also designed for posterior inference for the model parameters.

\emph{Heuristic}. Galbrun et al.~\cite{Galbrun_DMKD2014} proposed a generic greedy algorithm in finding overlapping communities for both labelled (on vertices) and non-labelled graphs, by maximising the edge density over all detected communities. They mapped the problem to the generalised maximum coverage problem~\cite{Cohen_IPL2008}, and obtain three different variants for their greedy algorithm by changing the greedy step. Li et al.~\cite{Li_WWW2015} on the other hand addressed overlapping CD problems in large networks based on a local spectra method where its running time relies solely on the size of the communities. The method shows high effectiveness and efficiency, given that it does not involve a huge number of computation on singular vectors, which enhance the weakness of traditional spectral clustering methods.

\emph{Non-negative matrix factorisation}. A framework that can be used in solving large scale CD problems was developed by Yang and Leskovec~\cite{Yang_WSDM2013}. They enhanced NMF in two aspects, aiming to estimate non-negative latent vector of each vertex and the likelihood formulation uses a near constant computation time. The running time of their framework is almost linear, and it discovers both overlapping and non-overlapping communities in a more accurate manner. Wang et al.~\cite{Wang_DMKD2011} proposed three approaches (asymmetric, symmetric and joint) based on NMF to detect hidden communities on various types of networks. They investigated the convergence and correctness of these algorithms, and showed their effectiveness using real-world networks.

We now discuss two learning-based models that have been established in addressing overlapping CD problems. The first model is used to discover both overlapping and non-overlapping communities in multi-view graphs whereas the second is developed based on GNNs.

A collection of graphs with multiple views where each view has the same vertex set but distinct edge set is defined as a \emph{multi-view} or \emph{multi-aspect} graph $ G $~\cite{Gujral_SDM2018}. As for ordinary graphs, each view in $ G $ is associated with an adjacency matrix. Hence, if $ G $ contains $ s $ views and $ n $ vertices, then it has $ s $ $ n \times n $ adjacency matrices.
\begin{problem}[\cite{Gujral_SDM2018}]\label{prob:Gujral2018}
	Take a multi-view graph $ G $ that has a small portion of vertex labels to $ c $ communities, assign all vertices of $ G $ to at least one of the $ c $ communities.
\end{problem}

In addressing Problem~\ref{prob:Gujral2018}, Gujral and Papalexakis~\cite{Gujral_SDM2018} suggested a semi-supervised framework to detect overlapping and non-overlapping communities. They called it as SMACD as the acronym of semi-supervised multi-aspect community detection.
%A parameter tuning algorithm that does not depend on the partial vertex labels was developed, so that SMACD is parameter-free.
%Let $ r $ be the number of communities.
The method first determines the decomposition of $ r-1 $ ($ r $, respectively) components for non-overlapping (overlapping, respectively) communities. For non-overlapping communities, each vertex is then assigned to one community that contains the maximum elements. In designing a new algorithm based on the coupled matrix-tensor factorisation~\cite{Acar_arxiv2011} model, they included two constraints, i.e.,non-negativity and latent sparsity, in dealing with CD.

Since the objective function of the latent sparsity constraint is hardly to be optimised and non-convex, the sparsity regulariser penalty $ \lambda $ is adapted. By convention, the trial and error approach is used to determine the parameter $ \lambda $, which will possibly affect the final outcome. To overcome this, the authors introduced an automated selection model that will automatically select an acceptable $ \lambda $, based on the fact that there exists a connection between $ \lambda $ and the levels of sparsity.

The SMACD model is evaluated by using eight real (two of them contain overlapping communities) and two synthetic datasets~\cite{Papalexakis_ICIFS2013}. Eight baselines~\cite{Papalexakis_ICIFS2013,Gligorijevic_ICPR2016,Nie_IJCAI2016} including GraphFuse (a special case of SMACD without supervision) are used for comparison purposes. Their code and the synthetic data are available online at \href{http://www.cs.ucr.edu/~egujr001/ucr/madlab/src/SHOCD.zip}{\texttt{http://www.cs.ucr.edu/{\small $ \sim $}egujr001}}.
%and all the data used for comparison purposes are also available online.
They evaluated the behaviour of each method according to three measures, they are purity, adjacent random index (ARI) and normalized mutual information (NMI). For most instances, their results indicate that SMACD outperforms other baselines. They also observed that the number of labels enhance the performance of the model for overlapping communities.

The authors also found that the degree of supervision is proportional to the achievement of the SMACD model. This rises a concern if a supervised model should be developed. For the automated selection model, its effectiveness is evaluated against a brute force approach. The value of $ \lambda $ that is obtained by both approaches in fact shows comparable accuracy. Hence, it may be worth to identify differences particularly the running time for these two methods.

In 2019, Shchur and G\"{u}nnemann~\cite{Shchur_DLG2019} introduced the neural overlapping community detection (abbreviated as NOCD) framework that combines GNNs with the Bernoulli-Poisson (BP) model~\cite{Yang_WSDM2013} to discover overlapping communities in undirected graphs. Note that BP model is an efficient graph generative model that can generate communities with various topologies, which is practical for overlapping CD.

Several benefits are discovered by employing GNN architectures in CD. For instance, compared to simpler models including MLP and free variable, outputs with better quality can be obtained by using GNN models given that they can produce affiliation matrices that are identical for adjacent vertices. In addition, the vertex features can be incorporated easily into the model. For unseen vertices during the training phase, it is also possible for them to predict their communities inductively.

As in~\cite{Gujral_SDM2018}, NMI is used as the evaluation metric because of its robustness and suitability in handling degenerate cases. A 2-layer GCN is used in their model for evaluation purposes. They observed that the performance of the model improves significantly by using batch normalisation after the first convolution layer, and applying $ L_{2} $ regularisation to all weight matrices. These two modifications set a distinction between their model and the conventional GCN.
%Note that Adam optimiser~\cite{Kingma_ICLR2015} is used for training purposes.
Their results show that the two variants of NOCD (with and without attributes) achieve the highest NMI in 90\% of the datasets. They observed that there exists a strong relationship in between the NMI of these two variants and the loss function. The model is also highly scalable even when the number of attributes is large.

\subsubsection*{Community Detection Using COVID-19 Dataset}

Due to the outbreak of COVID-19 that causes public health crisis around the world, Chaudhary and Singh developed unsupervised ML techniques in studying CD problems using COVID-19 dataset~\cite{Chaudhary_SNAM2021}. They analysed the trend and variations of infected cases in 187 countries, ranging from January to August 2020.
%The dataset is obtained from the website of John Hopkins University, which consists of the number of COVID-19 cases of 187 countries from January 22, 2020 to August 15, 2020.

To reduce the dimension of the dataset and identify the most useful variables, the authors employed the principal component analysis (PCA) so that most of the information will be retained. The covariance matrix for each of the 13 variables is generated to obtain its eigenvalues and eigenvectors. The principal component of the dataset is then the eigenvector with the largest eigenvalue. The eigenvalues that corresponds to the importance of the principal components are arranged in a descending order. Based on their dataset, the authors found that the first eight principal components are adequate for communities detection purposes given that a huge portion of the cumulative variance are derived from these variables.

To display communities from the heterogeneous elements, $ k $-means clustering (an unsupervised clustering framework) is applied on the reduced dataset. They determined the number $ k $ of clusters by employing the elbow method. Experimental results for the two variants (with and without the PCA) of $ k $-means clustering are analysed. Their results prove that $ k $-means clustering with the PCA achieves a better outcome in displaying communities.

\subsubsection*{Community Detection in Attributed Networks}	

CD problems for complex attributed networks aim to partition a network by considering not just the structure of the network, but also its vertex attributes where most vertices within a same group should possess the same attribute. They can be formulated as discrete optimisation problems by setting various objectives based on graph structures and vertex attributes. This problem has been addressed by using multi-objective evolutionary algorithms (MOEAs), which is known as one of the promising solvers for discrete optimisation problems. We first briefly discuss MOEAs that have been proposed in solving CD problems for attributed graphs, before extending it to the associated learning-based model.

Motivated by MOEAs that are widely used in handling conflicting objectives, Li et al.~\cite{Li_ITC2018} adapted this approach in solving CD problems for attributed graphs. They proposed a model MOEA-SA that considers two different objectives, based on both structural and attribute similarities. They employed a hybrid representation in order to fully utilise relationships among vertices, designed an operator that is used in guiding the evolutionary process, and introduced a correction strategy to handle improper solutions. By using similar concepts, Pizzuti and Socievole~\cite{Pizzuti_ITC2020} introduced MOGA-@Net that evaluates three objectives in optimising the structural properties and three vertex similarity measures in accessing feature homogeneity. A post local merge procedure is also proposed to reduce the number of communities to obtain high quality outputs.

Since it is often hard to perform a search in discrete problems due to limited information, Sun et al.~\cite{Sun_ITC2022} proposed a continuous encoding multi-objective evolutionary algorithm (CE-MOEA) by using a GNN encoding method, which transformed discrete problems to continuous problems. For the GNN encoding, each edge in a graph is first associated with a continuous value. For every vertex $ v $, let $ \mathbf{x}_{v} $ represents its associated continuous vector that contains all $ k $ edges incident to $ v $. A sigmoid function $ S(x) = \frac{1}{1+\mathrm{e}^{-x}} $ that can be thought as the sigmoid layer in the GNN is first applied on $ x \in \mathbf{x}_{v} $ to obtain $ \mathbf{h}_{v} $. The normalised exponential function (or softmax function) $ \sigma(h)_{i} = \frac{\mathrm{e}^{h_{i}}}{\sum_{j=1}^{k} \mathrm{e}^{h_{j}}} $ is then applied on $ h \in \mathbf{h}_{v} $ in the softmax layer, which gives a probability distribution $ \mathbf{p}_{v} $ where each $ p \in \mathbf{p}_{v} $ represents the probability of the associated vertex. The vertex is then selected in the next layer based on the arguments of the maxima, $ \mathbf{argmax} $.

To evaluate the network structure of a graph, the modularity~\cite{Newman_PRE2004} score is adapted as the objective function. To evaluate the vertex attribute similarity, the authors revised the two objective functions proposed in~\cite{Li_ITC2018} in handling both single- and multi-attribute networks. Their algorithm was developed based on a non-dominated sorting-based MOEA, namely NSGA-II~\cite{Deb_PRE2002}.

The strengths of their learning framework in solving CD problems for attributed graphs include: (1) the search strategy is more robust as the information of adjacency vertices are fully exploited at the softmax layer during the GNN encoding, (2) the encoding method can be used on undirected or directed graphs, regardless of whether they are attributed or non-attributed, and (3) by converting from discrete to continuous optimisation problems, any MOEA in solving continuous multi-objective optimisation problems can be employed, resulting in a smoother fitness landscape.

For evaluation purposes, they considered two single-attribute datasets and 13 multi-attribute datasets (five of them with ground truth), and evaluate the performance of their framework using three metrics: density, entropy and NMI. The first two metrics are meant for datasets without ground truth. The results demonstrate that their framework outperforms both EA-based and non-EA-based algorithms. The proposed objective functions are also appropriate compared to others.

\subsubsection*{Future Direction}

To extend the work for overlapping CD problems, one could study if the two variants of NOCD can be used in identifying a concrete relationship between attributes and structures of a given community. The inductive performance of the model could also be evaluated for further improvement~\cite{Shchur_DLG2019}. It is also natural to ask if NOCD can be modified in handling other complex networks such as multi-view graphs, as evidenced by its promising results.

To enhance the accuracy of the model proposed in~\cite{Chaudhary_SNAM2021}, more attributes could be considered in order to involve more variables in the dataset prior to the PCA. A thorough analysis in determining the $ k $ value should be elaborated, and relevant comparisons could be done to get a more informative outcome. Other reduction methods could be implemented to examine their suitability on CD problems using relevant datasets~\cite{Chaudhary_SNAM2021}.

Lastly, the GNN encoding method in studying CD problems for overlapping complex attributed networks could be extended. Li et al.~\cite{Li_ITC2018} also suggested to construct new GNN encoding methods in addressing other discrete optimisation problems. The objective functions of their model could also be studied and enhanced to further improve its performance.

\subsection{Community Search}\label{sec:community_search}

Suppose $ G $ is a graph and $ v \in V(G) $ is a query vertex, the objective of the \emph{community search} (CS) problem~\cite{Sozio_SIGKDD2010,Fang_VLDB2020} is to deduce a subgraph $ H \subseteq G $ that satisfies the cohesiveness and connectivity constraints, and $ v \in V(H) $. Classical cohesiveness metrics comprise $ k $-core~\cite{Seidman_SN1983}, $ k $-truss~\cite{Cohen_NSA2008}, $ k $-clique, $ k $-edge-connected component~\cite{Gibbons_Cambridge1985} as well as modularity approach~\cite{Kim_SIGMOD2022}.

The \emph{$ k $-core} of $ G $ is a maximal connected subgraph $ H_{1} $ of $ G $ such that $ \hbox{deg}(v) \ge k $ for each $ v \in H_{1} $. The \emph{$ k $-truss} of $ G $ is a maximal connected subgraph $ H_{2} $ of $ G $ such that every edge $ e \in E(H_{2}) $ belongs to more than $ k-1 $ triangles in $ H_{2} $. The \emph{$ k $-clique} of $ G $ is a connected subgraph $ H_{3} $ of $ G $ such that $ H_{3} $ is a complete graph of order $ k $. The \emph{$ k $-edge-connected component} of $ G $ is a connected subgraph $ H_{4} $ of $ G $ such that $ H_{4} $ remains connected if less than $ k $ edges are deleted from $ H_{4} $. CS can be applied in various real-world applications~\cite{Fang_VLDB2020} such as e-commerce, fraudulent group discovery and friend recommendation.

%The CS problem is closely related to the CD problem. For CS problems, the search is done in a given graph based on a set of query vertices but CD problems usually partition the whole graph to detect all communities. In CS problems, communities are defined according to the given parameters whereas global criteria are used in CD problems. Compared to CS problems, algorithms in solving CD problems are usually inefficient and hence not suitable for large graphs. 

Despite the fact that many CS algorithms~\cite{Sozio_SIGKDD2010,Huang_SIGMOD2014,Huang_arxiv2015,Fang_PVLDB2016,Fang_PVLDB2017,Huang_PVLDB2017,Jiang_ICDE2018,Fang_TKDE2019,Fang_VLDB2020,Fang_PVLDB2020} have been proposed, predefined subgraph patterns that lead to some ineffectiveness are typically selected in modelling communities that are rigid when applied to different scenarios. It is also challenging to determine a good $ k $ value for each subgraph cohesiveness metric. A detailed survey of this topic can be found in~\cite{Fang_VLDB2020}.

To enhance further the effectiveness, efficiency and flexibility of CS algorithms, learning-based approaches are recently proposed. We discuss four learning frameworks in addressing CS problems including interactive and attributed CS problems, as well as meta-learning in CS. Interactive CS can often be applied on social networks, attributed CS is a more challenging problem with additional features and meta-learning learns from multiple prior tasks in order to solve a new task by using a small amount of observed data.

\subsubsection*{Interactive Community Search}

There are some restrictions in solving CS problems in online social networks based on existing frameworks. A large portion of networks will usually be crawled prior to searching for the communities, which may involve a huge amount of unattractive data from users perspective. In addition, rules defined in determining the membership of each member are somehow less effective, given that various query vertices could be chosen.

To investigate CS problems for online networks, Gao et al.~\cite{Gao_PVLDB2021} developed an interactive CS algorithm that employs a GNN. It is abbreviated as ICS-GNN and it captures similarities between vertices concurrently in an online social network according to the content and structural features. In the model, a $ k $MG community, where $ k $ is an integer, is formed by modelling a community as a size-$ k $ subgraph that has the largest GNN scores. Their framework requires several rounds of CS, by incorporating users' feedback throughout the process. The process terminates when users are satisfied with the obtained communities.

Specifically, the breadth-first search algorithm is first applied to construct a subgraph that consists of a query vertex. The authors employed an exploration strategy to filter out irrelevant vertices during the construction of the subgraph. In order to deduce the probability for each vertex, they trained a model on the subgraph based on a GNN. The cross entropy is employed as the loss function and GNN scores of other vertices are deduced by the trained model.

To locate the $ k $MG community, some approximate algorithms are applied where the GNN score for each vertex is used to determine if a vertex belongs to a community. The trick is that vertices with higher GNN scores will be put into the community by swapping with vertices that have lower GNN scores, without breaking its connectivity. The swapping process stops when no more vertex pair can be swapped.

The GNN model that indicates if a vertex $ v $ should be in a community is trained by using labelled vertices provided by users. Hence, a ranking loss is proposed so that implicit feedback from users can be incorporated in order for them to figure out a correct label. For the specific case where $ v $ is equal to a boundary vertex, the authors used a greedy approach to identify the membership of $ v $ according to the shortest path in between $ v $ and any vertex in the community, by utilising the global relative benefit.

The efficiency and effectiveness of their model are analysed using various networks including an online social network namely Sina Weibo. The authors noticed that the proposed model expresses both content and structural features effectively, and the ranking loss is also proved to be effective for interactive CS. In addition, communities that are produced using the global relative benefit fit into various distributions.

\subsubsection*{Attributed Community Search}

The attributed community search (ACS) problem is another CS problem that have been studied extensively. In this problem, every vertex in a community is assigned with a homogeneous attributed value. Its objective is to determine a subgraph that contains all vertices that possess similar attribute to a given query vertex. Majority of the existing works such as ACQ~\cite{Fang_PVLDB2016} and ATC~\cite{Huang_PVLDB2017} solved this problem without considering the relationship between the query vertex and attribute function. They usually utilise a two-stage framework by first determining a dense community based on a query, and then contract it based on a given attribute function.

To address these limitations, Jiang et al.~\cite{Jiang_PVLDB2022} proposed a supervised framework that encodes various information from both graphs and queries, in solving both CS and ACS problems (with appropriate extensions). The framework is known as query-driven graph neural networks (QD-GNN). Three encoders and one feature fusion are designed in the learning model. This framework also extends ICS-GNN~\cite{Gao_PVLDB2021} in dealing with attributed graphs.

To obtain a robust model, the graph encoder in QD-GNN encodes the information of a query-independent graph where both vertex attribute features and graph structures are integrated. The query encoder on the other hand provides an interface for query vertices, and models their structure information. The feature fusion operator combines the embedding results generated by both graph and query encoders to obtain the final result of QD-GNN, by utilising the local query information and global graph knowledge.

By incorporating query attributes into QD-GNN, a new model AQD-GNN that can be used to solve the ACS problem is obtained by including one extra component attribute encoder and a revised fusion operator. The attribute encoder is designed mainly for query attributes, where a bipartite graph between the vertex and attribute sets is created to obtain the attribute vertex embedding, by using two different propagations. In the revised fusion operator, all three embeddings (derived from graph, query and attribute encoders) are combined in acquiring the final result of the problem.

The authors evaluated the behaviour of the model using 15 attributed graphs\footnote{\href{https://linqs.soe.ucsc.edu/data}{https://linqs.soe.ucsc.edu/data}} based on $ F_{1} $ score. The model QD-GNN and its variant AQD-GNN are compared with attributed~\cite{Fang_PVLDB2016,Huang_PVLDB2017} and non-attributed~\cite{Chang_SIGMOD2015,Huang_arxiv2015} CS algorithms, where their performances are in fact comparable. This may be because of the role of the fusion operation in transmitting graph information and query vertices to either the query or attribute encoder. To support interactive attributed CS, the authors replaced the GNN model in ICS-GNN~\cite{Gao_PVLDB2021} by QD-GNN and AQD-GNN. The experimental results validate the effectiveness and efficiency of the two new models. They also noticed the three encoders and feature fusion boost the effectiveness of the model drastically.

\subsubsection*{Contrastive Learning in Community Search}

Even though both ICS-GNN~\cite{Gao_PVLDB2021} and QD-GNN~\cite{Jiang_PVLDB2022} achieve promising results in solving CS problems, their performance is however restricted given that they require many labels during the training process. To overcome these limitations and enhance further the effectiveness and efficiency of learning models, Li et al.~\cite{Li_ICDE2023} developed a semi-supervised framework for undirected graphs, by proposing a partition model that requires just a small number of labels. The authors named the model as COCLEP and it is designed according to a deep learning technique namely contrastive learning.
%They claimed that this is the first time contrastive learning is used in CS problems.

To encode vertices of an undirected graph into low-dimensional vectors using the framework, the authors employed an attention-driven GNN that focuses on fusing the query vertex information and the graph structure. They introduced a label-aware contrastive learner to encapsulate information in the community, and hypergraphs are used and encoded by a hypergraph neural network (HGNN) to capture the structure information of neighbouring vertices. The HGNN helps to improve the weaknesses of existing augmentation techniques that fail to consolidate high-order graph structures and retain relationships between vertices. The graph is then partitioned based on its minimum cut, and a search is performed in each subgraph so that the community result can be obtained by merging them together. They showed that the partition is beneficial in reducing both the memory and time of the training.

They evaluated their framework using eight datasets that have ground-truth communities. Comparisons are made using both traditional and learning-based CS algorithms based on three different metrics. The experimental results demonstrate that COCLEP and its variant COCLE (a non-partition model) outperform their competitors in most instances. They also give promising results in terms of the training and query efficiencies.

\subsubsection*{Meta-Learning in Community Search}

Traditional CS algorithms are often designed based on predefined subgraph patterns, resulting in ineffective solutions particularly for real datasets. It is rather challenging to develop a universal framework due to a wide range of topological graph structures and query vertices. Learning-based methods~\cite{Gao_PVLDB2021,Jiang_PVLDB2022} based on ground truth communities have recently attracted a lot of attentions in enhancing these inadequacies. They are also expected to give a good generalisation on unknown vertex relationships.

One of the main considerations in designing ML methods is about the availability of samples. In dealing with the situation where only a limited number of samples are available, meta-learning plays an important role. Meta-learning aims to design learning-based algorithms based on prior learning so that they can be used in solving a new problem, by using a small amount of samples. A survey on meta-learning methods and applications with GNNs can be found in~\cite{Mandal_SIGKDD2022}.

Fang et al.~\cite{Fang_arxiv2022} proposed a metric-based meta-learning model that extracts information from small data. The proposed model, namely conditional graph neural process (CGNP), is an extension of the ordinary neural model, conditional neural process~\cite{Garnelo_ICML2018} (CNP), which generates vertex embeddings based on small training data. They studied three schemes by considering connected graphs with shared and disjoint communities, and disconnected graphs with disjoint communities.

There are three main components in CGNP, which includes a commutative pooling operation, GNN encoder and decoder. The encoder is a $ k $-layer GNN that maps a graph, a query vertex and the associated ground truth to a vertex embedding matrix. The embedding matrices for all query vertices are then combined into one representation based on three commutative operations that are self-attention~\cite{Vasmani_NIPS2017}, sum and average. Based on the combined representation, the decoder then predicts the membership of a new query vertex. Note that three inner product decoders are designed, they are a simple decoder, a multilayer perception decoder and a GNN decoder. The CGNP model is trained by using stochastic gradient descent, in optimising the negative log-likelihood.

Out of three well known GNN layers (GAT~\cite{Velickovic_ICLR2018}, GCN~\cite{Kipf_ICLR2017} and GraphSAGE~\cite{Hamilton_NIPS2017}), GAT is chosen because of its superior performance. They used accuracy, recall, precision and F1-score to evaluate performances of all variants of CGNP with other baselines. The experimental results demonstrate that CGNP and its variants achieve much better performances than traditional, supervised and meta-learning methods in all datasets, for both connected and disconnected graphs. Results for disconnected graphs and the performance of the model when it is trained using a small number of ground truth also validate the effectiveness of CGNP where it can learn from just a small amount of data. The proposed framework shows the best efficiency among learning-based methods, and it can handle graphs up to 10,000 vertices.

\subsubsection*{Future Direction}

Given the nature of ICS-GNN, one could investigate if active learning can be employed in addressing interactive CS. Other optimisation methods could also be used to enhance the labelling strategy, as well as to locate the communities. Since the model needs to be retrained for each query, a better training process should also be designed to reduce its computation cost and time.

For the ACS problem, it may be worth to enhance the model QD-GNN by limiting the number of training queries, particularly for large graphs, to improve further its efficiency. The framework could also be revised so that it works for other attributed graphs, e.g., location-based attributed graphs.

Given that meta-learning gives a promising result in solving CS problems, we suggest to extend this concept to CD problems. Since metric-based approaches only fit for classification tasks, one could consider to employ other meta-learning approaches including black-box adaption and optimisation-based meta-learning~\cite{Finn_ICML2017,Ravi_ICLR2017} in order to extend the framework to a wider range of graph problems.

\section{Learning-Based Methods for Classic Graph Problems}\label{sec:classic_graph_problems}

Knowing that there exist many learning frameworks that have been proposed in solving graph related problems in recent years, and to demonstrate the importance of this research direction, we now give a brief overview on learning frameworks that have been developed for six classic graph problems, mainly on combinatorial optimisation and NP-complete problems~\cite{Karp_1972,Garey_FREEMAN1979}. Relevant surveys can also be found in~\cite{Bengio_EJOR2021,Cappart_arxiv2021,Mazyavkina_COR2021,Peng_DSE2021,Karimi-Mamaghan_EJOR2022,Yow_arxiv2022}.

\emph{Maximum weight clique}. Sun et al.~\cite{Sun_TPAMI2021} introduced MLPR in solving the maximum weight clique (MWC) problem by shrinking the original problem using a greedy approach. In MLPR, supervised learning is employed and the training set is solved using an exact algorithm before the solution of a given hard problem is predicted by the model.

\emph{Minimum weight independent dominating set problem}. Wang et al.~\cite{Wang_IS2020} combined an RL-based repair procedure with a local search algorithm to tackle this problem. A method based on RL was first developed to enhance the local search procedure, and a repair procedure was developed so that the best solution can be obtained.

\emph{Maximum cut}. Dai et al.~\cite{Dai_NIPS2017} developed a framework based on graph embedding and RL, known as S2V-DQN in addressing the maximum cut problem. (The same framework can be modified to solve the travelling salesman and minimum vertex cover problems.). They presented a greedy framework and introduced a deep learning method structure2vec to parameterise evaluation functions over graphs. Both fitted $ Q $-iteration and $ n $-step Q-learning~\cite{Sutton_2018} are also merged to learn the algorithm.

\emph{Travelling salesman problem}. In solving the ordinary travelling saleman problem (TSP), Deudon et al.~\cite{Deudon_CPAIOR2018} introduced a model that relies on RL and the well known 2-opt heuristic to enhance the current best solution. As in~\cite{Bello_arxiv2016}, given the current solution, the policy gradient and neural networks are used to learn the best action. Another learning-based method in solving TSP is S2V-DQN~\cite{Dai_NIPS2017} that was mentioned in the previous paragraph, which focuses on RL and graph embedding. Prates et al.~\cite{Prates_AAAI2019} showed that GNNs with limited supervision can be used to tackle the decision TSP. Their model can predict this decision problem within a small deviation from the optimal tour cost that is computed by using the Concorde TSP solver~\cite{Applegate_CONCORDE2006}.

\emph{Graph colouring problem}. Zhou et al.~\cite{Zhou_COR2014} proposed to solve the graph colouring problem (GCP) by using a learning-based exact algorithm CDCL that uses propositional logic. Due to the complicated relationships between implicit constraints for non-adjacent vertices, the CDCL is enhanced so that three procedures are used in discovering these constraints in order to reduce the search space. Zhou et al.~\cite{Zhou_ESA2016} introduced a general framework RLS by combining both RL and local search techniques in solving the GCP. A probability smoothing technique based on forgetting mechanisms~\cite{Hutter_CP2002} is also developed to identify old decisions in the current search, which ultimately shows a faster convergence rate. Motivated by~\cite{Prates_AAAI2019}, Lemos et al.~\cite{Lemos_ICTAI2019} invented a basic model based on GNN in dealing with the decision GCP. They trained the model based on the stochastic gradient descent algorithm that is implemented via TensorFlow’s Adam optimiser~\cite{Kingma_ICLR2015}.

\emph{Minimum vertex cover}. Mousavian et al.~\cite{Mousavian_ICEE2014} proposed an algorithm that combines both cellular and learning automata in solving the minimum vertex cover (MVC) problem. They first represented every vertex of a graph by a cell that comes with a learning automaton, and the cell is penalised or rewarded according to some predefined criteria, until a given threshold is obtained. In solving the MVC problem using S2V-DQN, Dai et al.~\cite{Dai_NIPS2017} generated Barab\'{a}si-Albert~\cite{Albert_RMP2002} and Erd\H{o}s-Renyi~\cite{Erdos_1960} networks for experimental purposes. The approximation ratio of S2V-DQN is nearly optimal when comparing to other methods. The model also discovers a new algorithm for the MVC problem.

\section{Conclusions and Future Work}\label{sec:conclusion}

In this survey, we review 20 learning frameworks that have been proposed in addressing five well known subgraph problems, they are subgraph isomorphism (counting and matching), maximum common subgraph, community detection and community search problems. For each problem, we point out its traditional algorithms as well as ML models that have been designed. We discuss the structures of each model and analyse their performances when comparing with relevant baselines. The links of the datasets and source codes that are available online are also provided. We then suggest some potential problems and challenges for future research in the relevant section, so that more learning-based methods can be explored and extended to other graph problems.

For the sake of completeness and since learning models in this area are relatively limited, we now suggest some general questions related to this work, so that more learning frameworks can be designed in tackling relevant graph problems.

The \emph{density} of a subgraph $ H $ of a given graph is defined as $ d(H) = \frac{\abs{E(H)}}{\abs{V(H)}} $. The densest subgraph problem~\cite{Goldberg_TR1984,Asahiro_JOA2000,Asahiro_DAM2002,Khuller_ICALP2009,Ma_TODS2021} aims to search for a subgraph that has the maximum density. Many algorithms have been developed and there exist many variations for the densest subgraph problem. One of them is known as the \emph{densest $ k $-subgraph} problem that generalises the clique problem, where it determines the maximum densest subgraph that contains exactly $ k $ vertices in a given graph.
%The densest subgraph problem~\cite{Goldberg_TR1984,Asahiro_JOA2000,Asahiro_DAM2002,Khuller_ICALP2009,Ma_TODS2021} is one of the subgraph problems that is well studied over the years, in which the objective is to search for a subgraph that has the maximum density.
%Many algorithms have been developed in addressing this problem, for both undirected and directed graphs. There exist many variations for the densest subgraph problem.
To the best of our knowledge, there is currently no known learning model in solving the densest subgraph problem. Since there is a ML method~\cite{Sun_TPAMI2021} in addressing the MWC problem, and there is a relationship between them, it is hence natural to suggest:
\begin{problem}
	Can a learning framework be developed to address the densest subgraph problem, or any of its variants?
\end{problem}

Most of the methods~\cite{Sozio_SIGKDD2010,Huang_arxiv2015,Fang_PVLDB2016,Huang_PVLDB2017} in solving CS problems consider undirected static graphs. There are many real-world applications that are modelled by using directed graphs~\cite{Bondy_MACMILLAN1976,Bang-Jensen_SPRINGER2008,Yow_arxiv2018,Yow_GC2021}, to encode extra information based on the direction of each edge. Dynamic graphs~\cite{Kazemi_JMLR2020} on the other hand allow transactions or interactions over time when a system is modelled by using a graph. One such illustration is Twitter where vertices and edges will be created and removed regularly. These imply that directed and dynamic graphs in fact play a significant role, and hence more attention should be given on them. Knowing that directed graphs and dynamic graphs are somehow more restricted, existing work~\cite{Jiang_ICDE2018,Fang_TKDE2019} in dealing with these types of graphs in CS problems are relatively little. Therefore, learning-based approaches could probably be adapted.
\begin{problem}
	Can learning methods be developed in solving community search problems for other graph types particularly for directed and dynamic graphs, so that a more promising solution could be obtained?
\end{problem}

To solve a graph problem, different learning strategies could be proposed in ML frameworks.
\begin{problem}
	Can learning strategies be characterised to identify the most suitable approach that should be employed in tackling certain class of graph problems?
\end{problem}

Recall that the semi-supervised framework COCLEP has been shown to be useful in solving CS problems as it requires lesser labelled data, and the other three learning frameworks (see Section~\ref{sec:community_search}) in dealing with CS problems employ supervised learning. Since supervised frameworks rely heavily on the availability and quality of datasets, one may explore if non-supervised learning strategies are in fact more appropriate in this task. 
\begin{problem}
	Does the overall performances of learning frameworks in solving community search problems always be strengthened, by employing non-supervised ML strategies?
	%Can community search problems be solved by employing different ML strategies, rather than just supervised learning? If so, can the overall performances of the current frameworks be strengthened?
\end{problem}

%\paragraph{Paragraph headings} Use paragraph headings as needed.

%\begin{figure}
%  \includegraphics{example.eps}
%\end{figure}
%\begin{figure*}
%  \includegraphics[width=0.75\textwidth]{example.eps}
%\end{figure*}
%\begin{table}
%\begin{tabular}{lll}
%\hline\noalign{\smallskip}
%\end{tabular}
%\end{table}

\begin{acknowledgements}
	The first author is most grateful to Universiti Putra Malaysia for granting him Leave of Absence in completing this work in Nanyang Technological University, and Singapore National Academy of Science for appointing him as an SASEA Fellow through a grant supported by National Research Foundation Singapore. Siqiang Luo is supported by Singapore AcRF Tier 1 (RG18/21), and in part by AcRF Tier 1 Seed Funding (RS05/21) and NTU startup grant.
\end{acknowledgements}

\bibliographystyle{spbasic}      % basic style, author-year citations
\bibliography{references}   % name your BibTeX data base

\end{document}